\documentclass[twocolumn,10pt]{IEEEtran}

\usepackage{graphicx}
\usepackage{color}
\usepackage{epstopdf}
\graphicspath{{./fig/}}
\usepackage{amsmath,amsthm}
\usepackage {amssymb}
\usepackage{array}
\usepackage{url}
\usepackage{cite}
\usepackage{hyperref}
\usepackage{balance}

\newcommand{\tb}[0]{\color{black}}
\usepackage[colorinlistoftodos]{todonotes}

\DeclareMathOperator*{\argmin}{arg\,min}
\newtheorem{example}{Example}

\usepackage{tikz}
\usepackage{tikzscale}

  {\list{}{\leftmargin=0.16in\rightmargin=0.16in}\item[]}%
  {\endlist}

\DeclareGraphicsExtensions{.eps,.pdf,.png,.jpg}

\usetikzlibrary{arrows,arrows.meta}
\usetikzlibrary{calc}
\usetikzlibrary{fit}
\usetikzlibrary{shapes,shapes.arrows}
\usetikzlibrary{cd}

\tikzstyle{sum} = [draw, thick, inner sep=0mm, circle, minimum size=1mm]
\tikzstyle{gainright} = [draw, thick, isosceles triangle, minimum height = 3mm, isosceles triangle apex angle=60]
\tikzstyle{gainleft} = [draw, thick, isosceles triangle, minimum height = 3mm, inner sep = 0mm, isosceles triangle apex angle=60, shape border rotate=180]
\tikzstyle{gainup} = [draw, thick, isosceles triangle, minimum height = 8mm, inner sep = 0mm, isosceles triangle apex angle=60, shape border rotate=90]
\tikzstyle{none} = [draw=none]
\tikzstyle{connector} = [->,thick]
\tikzstyle{connectordirect} = [thick]
\tikzstyle{line} = [thick]
\tikzstyle{block} = [draw, rectangle, thick,minimum height=1em,minimum width=1em]
\tikzstyle{dashedblock} = [draw, rectangle, thick,dashed,minimum height=1em,minimum width=1em]
\tikzstyle{smallsum} = [draw,circle,inner sep=0mm,minimum size=3mm]
\tikzstyle{branch} = [draw,circle,inner sep=0.5mm,fill=black]

\tikzset{
  saturation block/.style={
    draw, thick,
    path picture={
      \pgfpointdiff{\pgfpointanchor{path picture bounding box}{north east}}
        {\pgfpointanchor{path picture bounding box}{south west}}
      \pgfgetlastxy\x\y
      \tikzset{x=\x*.4, y=\y*.4}
      \draw (-.9,0) -- (.9,0) (0,-.9) -- (0,.9);
      \draw (-.9,-.6) -- (-.6,-.6) -- (.6,.6) -- (.9,.6);
      \node[text width=.1cm] at (-.35,.55) {\scriptsize $P_\calU$};
    }
  }
}

\newcommand\oprocendsymbol{\hbox{$\square$}} \newcommand\oprocend{\relax\ifmmode\else\unskip\hfill\fi\oprocendsymbol}

\begin{document}

\title{
Towards a Systems Theory of Algorithms
}

\author{Florian D\"orfler\textsuperscript{*},
Zhiyu He\textsuperscript{*,$\dagger$}, Giuseppe Belgioioso\textsuperscript{*}, Saverio Bolognani\textsuperscript{*}, John Lygeros\textsuperscript{*}, \& Michael Muehlebach\textsuperscript{$\dagger$}
\thanks{* Automatic Control Laboratory, ETH Z\"urich, Switzerland. Email: \{dorfler,zhiyhe,gbelgioioso,bsaverio,jlygeros\}@ethz.ch. 
$\dagger$ Max Plack Institute for Intelligent Systems, T\"ubingen, Germany. Email: michael.muehlebach@tuebingen.mpg.de. 
This research is supported by the Max Planck ETH Center for Learning Systems and the Swiss National Science Foundation under the NCCR Automation.}
}

\maketitle

\begin{abstract} 
Traditionally, numerical algorithms are seen as isolated pieces of code confined to an {\em in silico} existence. However, this perspective is not appropriate for many modern computational approaches in control, learning, or optimization, wherein {\em in vivo} algorithms interact with their environment. Examples of such {\em open {\tb algorithms}} include various real-time optimization-based control strategies, reinforcement learning,  decision-making architectures, online optimization, and many more. Further,  even {\em closed} algorithms in learning or optimization are increasingly abstracted in block diagrams with interacting dynamic modules and pipelines.
In this opinion paper, we state our vision on a to-be-cultivated {\em systems theory of algorithms} and argue in favor of viewing algorithms as open dynamical systems interacting  with other algorithms, physical systems, humans, or databases. Remarkably, the manifold tools developed under the umbrella of systems theory {\tb are well suited for addressing a range\,of\,challenges in the algorithmic domain.}
We survey various instances where the principles of algorithmic systems theory are being developed and outline pertinent modeling, analysis, and design challenges.%
\end{abstract}


\section{Vision: Systems Theory of Algorithms} 
\label{sec:vision}

In the realm of control systems research, our traditional focus has predominantly revolved around dynamical systems, leveraging the underlying principles of physics, chemistry, biology, and so on to devise methods for modeling, analyzing, and controlling them. However, a notable shift has occurred due to the rise of computing techniques throughout the control stack, ranging from conventional microcontrollers with limited capabilities to advanced embedded real-time optimization and intricate networked control systems. 
From a systems theory perspective, the computing infrastructure has, to a large extent, been treated as a somewhat opaque entity. This perspective is gradually waning, prompting a necessary reevaluation, wherein algorithms themselves assume a pivotal role akin to open dynamical systems that are interconnected with other computing (or non-computing) systems in real time. 

This paradigm shift is {\tb neither entirely novel nor unique} to control systems. It {\tb has been recognized in specific contexts and is} also evident in optimization, machine learning, and other disciplines which increasingly conceptualize algorithms and computing pipelines as (interconnected) dynamical systems, reflecting a broader change in perspective.
We believe that the tools we have honed for analyzing and controlling dynamical systems 
can shed light on this emerging paradigm shift and help to navigate it. In fact, they have already proved themselves in many algorithmic challenges.
In this paper, we state our vision on a to-be-cultivated {\em systems theory of algorithms}. 

\subsection{Perspectives on Algorithms: Code or Dynamical Systems?} 
\label{subsec:perspectives}

We demarcate  two perspectives on algorithms that reflect some central watersheds in systems theory, such as the distinctions of closed versus open systems (i.e., that are either isolated or interacting with their environment), 
offline versus online decision making, or
black-box versus structured systems.
{\tb The two perspectives are intentionally overstated and caricatured to underscore our arguments more effectively.
In fact, many algorithms fall into the shades of grey between these extremes.}

{\tb The traditional viewpoint is that algorithms are merely {\bf  a piece of code} 
featuring some of the following characteristics:%
\begin{itemize}

\item closed, i.e., its evolution is autonomous in the dynamical systems sense: it is  based only on initialization, random seeds, and does not interact with the outside world; 

\item executed with batch data, not necessarily obeying causality, and operating offline or on a time scale distinct from how data from the ambient world is collected (e.g., offline smoothing of a previously recorded time series);

\item perfect in its execution, i.e., running until termination  determined by internal criteria, such as a relative tolerance or a maximum number of iterations; and

\item monolithic, i.e., appearing to the outside as a black box without any relevant internal structure.

\end{itemize}}

This perspective is commonly adopted, {\tb for instance,} when an $\text{\tt argmin}$ operation is employed as a subroutine in a controller.
In general, the above characteristics can be found in almost all implementations of numerical optimization-based control, 
in  control synthesis based on control barrier functions, in system identification, adaptive and learning-based control techniques, classical supervised learning methods, and so on.

In our opinion, this {\tb obviously caricatured} perspective of an algorithm as a ``piece of code'' is not suitable for many modern computational approaches in control, learning, optimization, game theory, and so on, wherein algorithms often need to be reactive rather than isolated. 
Here, we advocate the perspective that {\tb many algorithms take the form of a {\bf discrete-time dynamical system}} with some or all of the following characteristics:%
\begin{figure}[t]
	\centering{
	\includegraphics[width=1\columnwidth]{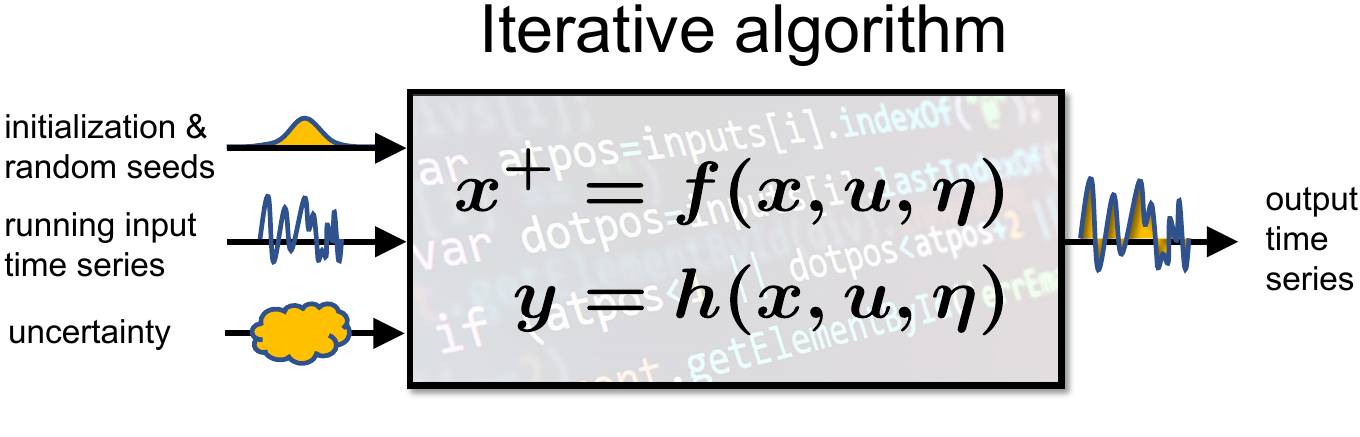}
		\vspace{-2em}
	\caption{
 We advocate modelling an algorithm as an open discrete-time dynamical system subject to inputs $u$, outputs $y$, an internal latent variable (state) $x$, and exogeneous  signal $\eta$ collecting different sources of uncertainty.}
	\label{Fig: open-alg}
	}
	\vspace{-1em}
\end{figure}
\begin{itemize}

\item open, i.e., endowed with inputs (e.g., real-time data or measurements) and outputs (e.g., {\tb convergence residuals}); 

\item  executed online with streaming data in a causal manner; 

\item subject to imperfections with (early) termination depending on external criteria (e.g., real-time requirements); and

\item possessing internal structure that can be leveraged when interacting with the algorithm (e.g., a monotone response of error residuals upon being queried by  inputs).

\end{itemize}
This alternative perspective, illustrated in Fig.~\ref{Fig: open-alg}, {\tb is valid for many numerical algorithms and} prompts several questions that can be studied with system-theoretic methods, including
\begin{itemize}

\item interconnection of algorithms, 
e.g., in machine learning pipelines, control stacks, parallel and distributed computing, or also in non-cooperative environments;

\item interconnection of algorithms with non-algorithmic (human or physical) systems, for example in recommender systems or real-time optimization-based control;

\item {\tb propagation of uncertainty, e.g., inexact algorithmic outputs (due to quantization or early stopping) become exogenous disturbances for downstream algorithms; or }

\item performance metrics beyond convergence, such as transient metrics (e.g., running cost or regret) or input/output gains (e.g., for disturbance amplification).

\end{itemize}
The reader can possibly already match different examples to these two paradigms and the different bullet items above. {\tb Our initial presentation of these paradigms is simplified and deliberately caricatured. Multiple real and in-depth case studies are presented in Section~\ref{sec:examples} showing the relevance of the  systems perspective and illustrate the above bullet items for certain classes of algorithms, old and new. In addition,} we present a simple and disarming example below, well-known to all control researchers and engineers, that can  reflect either the ``piece of code'' or the ``dynamical system'' paradigm.

\subsection{A Disarming Example Reflecting the Dichotomy}
%
Consider the problem of estimating an unknown quantity from a measurement residual (as in observer design, identification, adaptive control, and so on), which can be performed either by  least-squares minimization or with a Kalman filter \cite{kailath2000linear}. The two approaches reflect  the  characteristics of the ``piece of code'' and ``dynamical systems'' paradigms, respectively.

To be specific, consider the linear measurement equation
$y_{i}=C_{i}x + \eta_{i}$, where $x \in \mathbb R^{n}$ is the to-be-inferred parameter, and $i \in \{1,\dots, n\}$ indexes a sequence of $n$ measurements of $y_{i} \in \mathbb R$ and $C_{i} \in \mathbb R^{1 \times n}$ subject to independent zero-mean Gaussian noise $\eta_{i} \sim \mathcal N(0,\Sigma_{i})$ with variance $\Sigma_{i}>0$. 
In a Bayesian setting,  assume that a prior of $x \sim \mathcal N(0,\Pi)$ is given with positive definite $\Pi \in \mathbb R^{n \times n}$. Given batch data $(y_{i},C_{i})$, the {\em best} (maximum a posteriori probability) estimate of $x$ is obtained by minimizing the {\em regularized least-squares criterion}
\begin{equation}
\label{eq:LS}
	\hat x = \argmin\nolimits_{x}\, \sum\nolimits_{i=1}^{n} \left\| y_{i}- C_{i}x\right\|_{\Sigma_{i}^{-1}}^{2} + \left\| x  \right\|_{\Pi^{-1}}^{2}
	\,,
\end{equation}
the solution of which is given by the familiar regularized and weighted pseudo-inverse. This procedure reflects the perspective of the estimation algorithm to be a perfect and isolated piece of code running offline and with batch data.

Consider now a continuous data stream $(y_{i},C_{i})$, $i \in \mathbb Z_{\geq 0}$. In this setting, repeatedly solving the least-squares  problem \eqref{eq:LS} with ever increasing data size and ever larger matrices is impractical. Adopting a system-theoretic perspective, we model the data stream by means of the stochastic process
\begin{equation*}
x_{i} = x_{i-1} \,,\quad y_{i} = C_{i}x_{i} + \eta_{i}\,.
\end{equation*}
We can infer the best state estimate $\hat x_{i} \in \mathbb R^{n}$ and its covariance $P_{i} \in \mathbb R^{n\times n}$ at time $i$ by means of the Kalman filter
\begin{subequations}\label{eq: Kalman filter}
\begin{align}
	\hat x_{i} &= \hat x_{i-1} +  P_{i-1} C_{i}^{\top} \frac{y_{i}-C_{i}\hat x_{i-1}}{ \Sigma_{i} + C_{i} P_{i-1} C_{i}^{\top}}\,,
	\\
	 P_{i} &= P_{i-1} - \frac{P_{i-1}C_{i}^{\top}C_{i}P_{i-1}}{\Sigma_{i}+C_{i}P_{i-1}C_{i}^{\top}}\,,
\end{align}
\end{subequations}
initialized with $\hat x_{0} = 0$ and $P_{-1} = \Pi$.
\begin{figure}[t]
	\centering{
	\includegraphics[width=1\columnwidth]{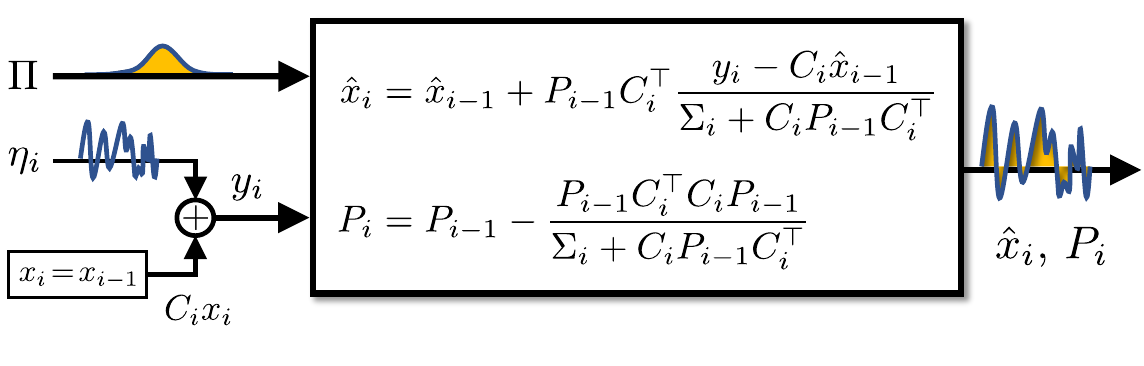}
		\vspace{-2em}
	\caption{
 The Kalman filtering approach \eqref{eq: Kalman filter} is an online algorithm to solve the least-squares problem \eqref{eq:LS} with streaming data.}
	\label{Fig: Kalman}
	}
		\vspace{-1em}
\end{figure}
The Kalman filter \eqref{eq: Kalman filter}, illustrated in Fig. \ref{Fig: Kalman}, is an online algorithm to solve the least-squares problem \eqref{eq:LS} with streaming data, also known as {\em recursive least squares}, and it  reflects the dynamical systems perspective on algorithms, i.e., it is online, causal, it is fed by streaming data, its structure is interpretable, it can be paired with downstream (e.g., control) algorithms, and so on.%

\subsection{Position} 
\label{subsec:statement}

In this opinion paper, 
we champion the view of algorithms as open dynamical systems, arguing that this paradigm is better suited for
the analysis and design of modern computational methods, computing pipelines, and their interaction with other computing systems, physical systems, or humans. 

Such thoughts -- advocating a systems theory perspective on algorithms -- have often been voiced from a didactic perspective \cite{choppella2021algodynamics} and for analysis and design of algorithms in optimization \cite{scherer2023optimization,lessard2022analysis,helmke2012optimization,brockett1991dynamical}, machine learning \cite{liu2019deep,LeeHe2020,li2020dynamical,Meyn2023}, and more general iterative algorithms \cite{bhaya2006control,hasan2012control,kashima2007system,chu2008linear,glover1973stability,stuart1998dynamical} in linear algebra, numerical integration, and so on. From a bird's eye view, the running theme in this literature is that fixed-point iterations --- arguably the workhorse of computing --- can be abstracted as a discrete-time dynamical system $x^{+} = f(x)$ (or, in continuous time, as a differential equation), typically with a structured right-hand side and amenable to a system-theoretic analysis. 

These prior works almost exclusively study algorithms as closed systems, whereas here we want to particularly advocate the viewpoint of an algorithm as an {\em open}  system interacting with its environment. An early and successful precursor in this line of thinking is the literature surrounding real-time iteration \cite{diehl2005real} analyzing Model Predictive Control (MPC) \cite{borrelli2017predictive} as the closed-loop interconnection of a (physical) dynamical system with a running Newton-type algorithm for solving the optimal control problem and demonstrating joint contraction of the physical and algorithmic dynamics. {\tb Further early touching points of control systems and algorithms are found in sampled-data and quantized control \cite{chen2012optimal,ishii2003quadratic}, {\tb extremum seeking \cite{ariyur2003real,scheinker2024100},}} error propagation of recursive least squares \cite{ljung1985error}, or the broader field of adaptive control \cite{annaswamy1989stable,aastromadaptive} seeking joint convergence of parameter estimation and a closed-loop system. 

On an abstract level, the field of control is concerned with the cyber-physical interconnection of controllers (implemented as digital ``algorithms'') with real-world processes.
Beyond control, {\em open} algorithms interacting with their environment (either other algorithms or non-algorithmic systems) are also common in the computing domain (parallel, distributed, or pipelines of algorithms), learning methods, and their applications; see Section~\ref{sec:examples} for a few contemporary examples.

In this opinion paper, we put forward {\bf three key messages} for the systems and control community:
\begin{itemize}

\item Some of our core strengths are the power of abstraction and deep understanding of feedback loops. We are trained to think in terms of block diagrams thereby imposing structure, modularizing, and taming complexity. Systems theory brings not only a different perspective but also a powerful set of tools to study and design  algorithms.

\item 
Modern algorithms are seldomly isolated pieces of code confined to an {\em in silico} existence, but they operate {\em in vivo}, i.e., in feedback with other algorithms, data bases, the physical world, or even humans. The concept of feedback is {\em the} very core of control, and that is why our community should embrace  topics in an ever-more digital, computational, data-rich, and algorithmic world.

\item 
While we can offer a new perspective, 
a direct application of our tools is insufficient for generating a lasting impact.
Therefore, we must expand our methodological base by integrating tools that are not broadly adopted in our field. These include 
computational tools (such as automatic differentiation and parallelization), analysis methods (such as (monotone) operator theory), statistical tools (such as sampling-based methods, generative modeling, or bootstrapping), and complexity theory, among others.

\end{itemize}

While the first two points encourage our involvement and urge us to prioritize algorithms as a central topic within the systems control community, the third point comes with a note of caution. It advises against unsolicited involvement and warns about the pitfalls of merely importing computational challenges into our field  for the sake and joy of mathematical exploration alone.
Instead, to create meaningful contributions, we need to adopt terminologies and techniques from the realm of computer science and engage in genuine collaborations.
%



\subsection{Outline} 
\label{subsec:outline}

In  Section \ref{sec:examples}, we present a suite of examples and success stories, where the dynamical systems perspective is appropriate and possibly the most effective approach to handle the problem's complexity. In Section~\ref{sec:challenges}, we present a list of grand challenges that can possibly be addressed  from the dynamical systems perspective. While these are certainly ambitious, they should guide and inform future research endeavors. Finally, Section \ref{sec:conc} concludes the paper with a brief discussion.

{\tb We close with two brief disclaimers: First, given the numerous touching points of system theory and algorithms, we do not intend to comprehensively survey the vast body of related work -- MPC itself is a universe next to many others -- but we selectively sample the literature. Second, our advocated position on a to-be-cultivated systems theory of algorithms is neither novel nor unique: many of the quoted references date back decades. Further, system-theoretic tools like Lyapunov, small gain, or passivity are widely used in the algorithmic domain, often under different names: for instance, the interplay between monotone operators and passivity properties is well-documented \cite{PavelDissipativity2022}. However, the system-theoretic perspective on algorithms is increasingly timely and shared, as shown next.} 


\section{When the System-Theoretic Perspective on Algorithms Proves Advantageous} 
\label{sec:examples}

In this section, we provide a set of compelling examples that demonstrate the application of systems theory in examining open and interconnected algorithms operating in  {\em in vivo} or {\em in silico}  environments, that is,   interconnected with physical systems or humans or confined to the computing world (though in setting where a system-theoretic viewpoint proves very useful). 
We selected three domains to highlight the convergence and synergistic interaction between systems theory and computational algorithms: the analysis, design, and interplay of optimization and learning algorithms; real-time algorithms in feedback loops; and decision-making architectures.

\subsection{Analysis, Design, and Interplay of Algorithms in Optimization and Learning}

The system-theoretic perspective on algorithms has already had a significant impact on the analysis and design of optimization and learning algorithms in the computing world. This influence is likely to become even more dominant, as emerging control and machine learning challenges -- such as personalized health care, autonomous supply chain management, or the development of recommender systems -- continue to grow in complexity. Furthermore, as algorithms become increasingly complicated and intertwined, their interplay calls for an effective means of characterizing and taming complexity.

The system-theoretic view has the potential to address these challenges by providing the following important features: 
\begin{itemize}
    \item \emph{Abstraction:} Abstracting algorithms as systems allows us to think in block diagrams. This breaks down complexity and enables the study of feedback, uncertainty, bias propagation, and dynamic interaction among blocks.

    \item \emph{Analysis and design:} Systems theory provides concepts, such as time-scale separation, system gains (e.g., $\mathcal{H}_2/\mathcal{H}_\infty$), integral quadratic constraints, and small-gain theorems, which provide insights into the stability and performance of algorithms, and are useful for design -- especially when algorithms interact with one another or with their environment; see also Section~\ref{eq:RTAFL}. 

    \item \emph{Going beyond convergence:} Traditional analysis often focuses on convergence certificates, but systems theory allows us to answer broader questions about disturbance rejection, regret, or architectural design, see also Section~\ref{subsec:architectures}. This shift in focus can be crucial for real-world applications, where convergence might be insufficient.
\end{itemize}

Mathematical optimization has been strongly influenced by systems theory ideas. An early illustration is the Arrow-Hurwitz method\cite{arrow1958studies} for computing saddle points, which was introduced as a dynamical system and analyzed with tools from calculus. Later, researchers such as R. Brockett, U. Helmke, and J. B. Moore \cite{brockett1991dynamical,helmke2012optimization} expanded the use of calculus and differential geometry to demonstrate that computational tasks, such as diagonalizing matrices or sorting lists, can be framed as isospectral flows, i.e., matrix differential equations preserving eigenvalues. Recent investigations into momentum-based optimization \cite{muehlebach2021optimization} have also highlighted connections to symplectic geometry, enabling a qualitative understanding of momentum and establishing convergence rates in a nonconvex setting. Moreover, the analysis and design of gradient-based algorithms {\tb and numerical integrators for ordinary differential equations \cite{glover1973stability,stuart1998dynamical}} are closely related to the Lur'{e} problem, i.e., the feedback interconnection of a linear dynamical system with a static nonlinearity \cite{khalil2002nonlinear}.
This relation led to many works deriving {\em i)} tight upper bounds on the convergence rate via integral quadratic constraints \cite{lessard2022analysis,scherer2023optimization}, parametric Lyapunov functions \cite{van2023tutorial}, or the Performance Estimation Problem approach\cite{goujaud2023fundamental,rubbens2023interpolation}, {\em ii)} lower bounds \cite{muehlebach2020continuous}, and {\em iii)} convergence rates for distributed optimization\cite{sundararajan2020analysis}.
It also inspired the principled design of gradient-based optimization algorithms by tuning algorithm parameters via system-theoretic tools\cite{van2017fastest}.

The following examples illustrate  {\tb this systems perspective} on prototypical optimization problems and {\tb (continuous- and discrete-time)} algorithms to solve them. 
In the first example, we present the role of feedback control in shaping the performance of a primal-dual optimization algorithm.

\begin{example}[Primal-dual algorithms as proportional-integral controllers]
	We discuss the control-theoretic interpretation of the role of the augmented Lagrangian\cite{wang2011control}. Consider the following equality-constrained problem
		\begin{equation}\label{eq:constrained_opt}
		\begin{split}
			\min\nolimits_{x} ~ f(x) \qquad \textup{s.t.} ~ Ax - b = 0,
		\end{split}
		\end{equation}
		where $x\in \mathbb{R}^n$, $A \in \mathbb{R}^{m\times n}$, $b \in \mathbb{R}^{m}$, and $f: \mathbb{R}^n \to \mathbb{R}$ is the objective. The augmented Lagrangian is $L(x,\lambda) = f(x) + \lambda^{\top}(Ax-b) + \tfrac\rho2 \|Ax-b\|^2$, where $\lambda \in \mathbb{R}^m$ is the multiplier, and $\rho>0$ is the penalty parameter. To solve problem~\eqref{eq:constrained_opt}, the classic saddle-point flow associated with $L(x,\lambda)$ is
		\begin{subequations}\label{eq:saddle_point_dynamics}
		\begin{align}
			\!\dot{x} &\!=\!-\nabla_x L(x,\lambda) = -\nabla f(x) - A^{\top}\lambda -\rho A^{\top}(Ax-b), \label{eq:saddle_point_primal} \\
			\!\dot{\lambda} &=\!+\nabla_{\lambda} L(x,\lambda) = Ax - b. \label{eq:saddle_point_dual}
		\end{align}
		\end{subequations}
		From a dynamics perspective, the vector field \eqref{eq:saddle_point_dynamics} decomposes into symplectic (rotational) dynamics induced by the constraint $Ax-b$ and dissipation induced by the gradient of the cost $f(x)$ aided by the quadratic {\tb penalty} $\tfrac\rho2 \|Ax-b\|^2$.
		From an optimization perspective, this quadratic {\tb penalty} regularizes $f(x)$ which is not necessarily strongly convex, thus facilitating convergence.
		Finally, from a control viewpoint, this quadratic {\tb penalty} yields a proportional control term in \eqref{eq:saddle_point_primal} that regulates the transient constraint violation, i.e., $Ax-b$. The dual ascent \eqref{eq:saddle_point_dual} further brings an integral part, thus contributing to constraint satisfaction.
		The block diagram in Fig.~\ref{fig:diagram_aug_lag} reveals the structure of the overall dynamics \eqref{eq:saddle_point_dynamics} as a proportional-integral (PI) control. This structure results in favorable properties, such as fast stabilization and elimination of steady-state error.

		\begin{figure}[!t]
			\centering
		    \resizebox{.7\columnwidth}{!}{
		    \begin{tikzpicture}
		      \matrix[ampersand replacement=\&, row sep=0.2cm, column sep=.9cm] {
		        \node[smallsum](sum1){}; \& \node[branch](britr) {}; \& \node[block](proportion) {$\rho$}; \& \node[smallsum](sum2){}; \& \\
		        \node[block](A) {$A$}; \& \& \node[block](integral) {$\int$}; \& \& \node[block](A_transpose) {$A^{\top}$}; \\
		        \& \node[branch](britr2) {}; \& \node[block](integral2) {$\int$}; \& \node[smallsum](sum3){}; \& \\
		        \& \& \node[block](grad) {$\nabla f$}; \& \& \\
		      };

		          \draw[connector] ([xshift=-.5cm]sum1.west)--(sum1.west) node[at start, left]{$b$} node[at end, above left] {$-$};
		          \draw[connector] (sum1.east)--(proportion.west);
		          \draw[connector] (britr.south)|-(integral.west);
		          \draw[connector] (proportion.east)--(sum2.west) node[at end, above left]{$+$};
		          \draw[connector] (integral.east)-|(sum2.south) node[at end, below right]{$+$} node[at start, above right]{$\lambda$};
		          \draw[connector] (sum2.east)-|(A_transpose.north);
		          \draw[connector] (A_transpose.south)|-(sum3.east) node[at end, above right]{$-$};
		          \draw[connector] (sum3.west)--(integral2.east);
		          \draw[connector] (britr2.south)|-(grad.west);
		          \draw[connector] (grad.east)-|(sum3.south) node[at end, below right]{$-$};
		          \draw[connector] (integral2.west)-|(A.south) node[at start, above left]{$x$};
		          \draw[connector] (A.north)--(sum1.south) node[at end, below left]{$+$};

		          \node[fit=(britr)(proportion)(integral)(sum2), draw, dashed, inner xsep= 2.5mm, inner ysep = 1mm, thick] (fitb1){};
    
				  \node[above right, inner sep = 1mm] at (fitb1.north west) {{\small PI control}};
		    \end{tikzpicture}
		    }
		    \vspace{-.5em}
		    \caption{The block diagram illustrates the saddle-point dynamics \eqref{eq:saddle_point_dynamics}.}
    		\label{fig:diagram_aug_lag}
		\vspace{-0em}
		\end{figure}
\end{example}

The second example showcases how a broad class of algorithms can be abstracted via a classical system-theoretic framework, namely as instances of the Lur'e problem \cite{khalil2002nonlinear}.

\begin{example}[Gradient-based algorithms as closed-loop systems]
    We illustrate how a system-theoretic perspective on gradient-based optimization algorithms facilitates analysis and design\cite{lessard2022analysis,van2017fastest}. Consider the  unconstrained problem
	\begin{equation}\label{eq:unconstrained_opt}
		\min\nolimits_x ~ f(x),
	\end{equation}
	where $f: \mathbb{R}^n \to \mathbb{R}$ is the objective. A multitude of gradient-based optimization algorithms for solving problem \eqref{eq:unconstrained_opt} are based on the following parameterized dynamics
	\begin{equation}\label{eq:gd_alg_family}
	\begin{split}
		\xi_{k+1} &= (1+\beta)\xi_k - \beta \xi_{k-1} - \alpha \nabla f(y_k), \\
		y_k &= (1+\gamma)\xi_k - \gamma \xi_{k-1}, \\
		x_k &= (1+\delta)\xi_k - \delta \xi_{k-1},
	\end{split}
	\end{equation}
	where $x_k \in \mathbb{R}^n$ is the primal variable, $\xi_k, y_k \in \mathbb{R}^n$ are memory states, and $\alpha>0$ and $\beta,\gamma,\delta \geq 0$ are design parameters\cite{van2017fastest}. For instance, $\beta = \gamma = \delta = 0$ recovers gradient descent, and  Nesterov's accelerated gradient method \cite{nesterov1983method} is obtained for $\beta = \gamma$ and $\delta = 0$. Figure~\ref{fig:diagram_gd_alg_family} illustrates the abstraction of algorithm \eqref{eq:gd_alg_family} into the feedback interconnection of a linear dynamical system and a static nonlinearity, the gradient $\nabla f$, which is a prototypical instance of the Lur'e problem. For the analysis integral quadratic constraints can be leveraged to characterize this nonlinearity and to derive sufficient parametric conditions that ensure a stable closed-loop system (i.e., a convergent algorithm) when the objective $f$ belongs to a certain class \cite{lessard2022analysis}. In terms of design, by analyzing the zeros and poles of the transfer functions of the blocks in Figure~\ref{fig:diagram_gd_alg_family}, a set of parameters for \eqref{eq:gd_alg_family} can be identified that leads to the so-called triple momentum method, which enjoys the fastest known convergence rate\cite{van2017fastest}.

	\begin{figure}
		\centering
		    \resizebox{.95\columnwidth}{!}{
		    \begin{tikzpicture}
		      \matrix[ampersand replacement=\&, row sep=0.3cm, column sep=.55cm] {
		        \node[none](edgetop1) {}; \& \node[block](algorithm) {$\begin{aligned}
		        \left[
		        \renewcommand \arraystretch{1.2}
		        \begin{array}{c}
		             \xi_{k+1} \\ \xi_k \\ \hline y_k \\ x_k
		        \end{array}
		        \right] =
		        \left[
		        \renewcommand \arraystretch{1.2}
		        \begin{array}{cc|c}
		            (1+\beta)I_n & -\beta I_n & -\alpha I_n \\
		            I_n & 0 & 0 \\
		            \hline
		            (1+\gamma)I_n & -\gamma I_n & 0 \\
		            (1+\delta)I_n & -\delta I_n & 0 \\
		        \end{array}
		        \right]
		        \cdot
		        \left[
		        \renewcommand \arraystretch{1.2}
		        \begin{array}{c}
		             \xi_k \\ \xi_{k-1} \\ \hline u_k \\
		        \end{array}
		        \right]
		        \end{aligned}$}; \& \node[none](edgetop2) {}; \\
		        \\
		         \& \node[block, minimum width = 3.5cm, minimum height = 0.9cm](feedback) {$u_k = \nabla f(y_k)$}; \& \node[none](edgetop4) {}; \\
		      };

		      \draw[connector] (edgetop1.west)--(algorithm.west); 
		      \draw[connectordirect] (feedback.west)-|(edgetop1.west) node[left, yshift = -1.1cm]{$u_k$};
		      \draw[connector] (edgetop4.east)--(feedback.east); 
		      \draw[connectordirect] (algorithm.east)-|(edgetop4.east) node[at end, right, yshift = 1.1cm]{$y_k$};
		    \end{tikzpicture}
		    }
		    \vspace{-.5em}
		    \caption{The block diagram illustrates the structure \eqref{eq:gd_alg_family} of gradient-based optimization algorithms, where $I_n$ denotes the identity matrix of size $n$.}
    		\label{fig:diagram_gd_alg_family}
		\vspace{-.5em}
	\end{figure}
\end{example}

Whereas the previous examples revealed a natural feedback control structure emerging in an isolated optimization algorithm confined to the computing world, we now turn to a feedforward viewpoint useful in the design of algorithms that are subjected to disturbance inputs from their environment.

\begin{example}[Feedforward for time-varying optimization]
    Feedforward {\tb methods} can be used for accurate tracking of the solutions to time-varying optimization problems of the form 
    \begin{equation}\label{eq:time_varying_opt}
        x^*(t) = \argmin\nolimits_{x\in \mathbb{R}^n}~ f(x;t),
    \end{equation}
    where the objective $f(x;t)$ is twice continuously differentiable and strongly convex in $x$, and it is revealed before the decision $x$ is committed at every time $t$. Such a time-varying objective can result from real-time decision-making with streaming data (e.g., disturbances or tracking signals) from the environment. 
  
    The standard Newton flow to solve problem~\eqref{eq:time_varying_opt} is 
    \begin{equation}\label{eq:newton_flow_TV}
        \dot{x}(t) = - \nabla_{xx}f(x;t)^{-1} \nabla_x f(x;t)\,,
    \end{equation}
    where $\nabla_{xx} f(x;t)$ is the Hessian matrix. Implementing this solution will result in a tracking error that depends on $\dot{x}^*(t)$, the rate of change of the optimizer. {\tb To achieve asymptotically accurate tracking, we can exploit a feedforward correction term that accounts for the variation $\dot{x}^*(t)$.} We know from the implicit function theorem that $\dot{x}^*(t) = -\nabla_{xx} f(x^*;t)^{-1}\nabla_{tx}f(x^*;t)$, where $\nabla_{tx}f$ is the derivative of the gradient $\nabla_x f$ with respect to time. 
    {\tb Since $x(t)$ instead of $x^*(t)$ is available, a running estimate of $\dot{x}^*(t)$ is $-\nabla_{xx} f(x;t)^{-1} \nabla_{tx} f(x;t)$. Thus, an algorithm that utilizes feedforward correction via such an estimate of $x^*(t)$ is}
    \begin{equation}\label{eq:pred_cor_alg}
	   \dot{x}(t) = -\nabla_{xx} f(x;t)^{-1} \bigl(\nabla_x f(x;t) + \nabla_{tx} f(x;t) \bigr) \,.
    \end{equation}
    With \eqref{eq:pred_cor_alg}, $x(t)$ asymptotically converges to $x^*(t)$ and tracks the solution of \eqref{eq:time_varying_opt} \cite{simonetto2020time,picallo2022sensitivity}.
		
\end{example}

Beyond the above success stories in optimization, systems theory has also triggered important developments in machine learning. Novikoff's proof of the convergence of the perceptron algorithm\cite{novikoff1962convergence}, that describes the perceptron as a dynamical system, can be seen as the birth of statistical learning theory. Similarly, the backpropagation algorithm is equivalent to the adjoint or co-state dynamics of a particular optimal control problem \cite{lecun1988theoretical}. Co-state dynamics also find application in neural ordinary differential equations, where the output corresponds to the solution of an ordinary differential equation with the right-hand side specified by a neural network. Furthermore, the core of diffusion-based generative modeling lies in a stochastic differential equation that links a prior (e.g., Gaussian) distribution with an unknown data distribution. 
{\tb This link elucidates that the Kullback--Leibler divergence\cite{cover2006elements} between the data distribution and the Gaussian distribution corresponds to the minimal total control cost\cite{tzen2019theoretical}, and it facilitates deriving sampling algorithms and the optimal drift \cite{tzen2019theoretical,berner2022optimal}.} 
{\tb Finally, the growing integration of state-space models into the architectures of foundational deep neural networks demonstrates the influence of systems theory \cite{alonso2024state}.}

\begin{figure}[!t]
        \centering
        \resizebox{\columnwidth}{!}{
        \begin{tikzpicture}
          \matrix[ampersand replacement=\&, row sep=.6cm, column sep=.7cm] {
            
            \node[branch](britr) {}; \& \node[block](preprocess) {$\begin{gathered}
                \text{Pre-} \\
                \text{processing}
            \end{gathered}$}; \& \node[block] (ssm) {
              $\begin{aligned}
                 x_{k+1} & = A_k x_k + B_k u_k \\
                  y_k & = C_k x_k + D_k u_k
                \end{aligned}$
            }; \& \node[block](postprocess){$\begin{gathered}
                \text{Post-} \\
                \text{processing}
            \end{gathered}$}; \\ 
             \& \& \node[block, minimum width = 2.5cm, minimum height = 0.9cm](activation) {Activation}; \&  \\
          };
    
          \draw[connector] ([xshift=-.5cm]britr.west)--(preprocess.west) node[at start, left]{$\tilde{u}_k$};
          \draw[connector] (preprocess.east)--(ssm.west) node[at end, above left] {$u_k$};
          \draw[connector] (ssm.east)--(postprocess.west) node[at start, above right] {$y_k$};
          \draw[connector] (britr.south)|-(activation.west);
          \draw[connector] (activation.east)-|(postprocess.south);
          \draw[connector] (postprocess.east)--([xshift=.5cm]postprocess.east) node[at end, right]{$\tilde{y}_k$};
    
        \end{tikzpicture}
        }
        \vspace{-1.5em}
        \caption{{\tb The block diagram adapted from \cite[Figure 1]{alonso2024state} shows the structured state-space sequence model. Stacks of such models form a foundation model.}}
        \label{fig:mamba-s6-block}
        \vspace{-1em}
\end{figure}
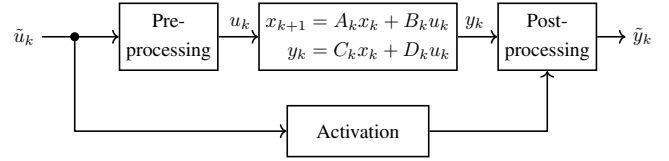

{\tb
\begin{example}[Structured state-space sequence model]
    Foundation models are widely adopted in modern machine learning practice. These are pre-trained on massive data and then fine-tuned for downstream tasks. To date, these models have predominantly relied on the {transformer architecture} and the attention mechanism\cite{vaswani2017attention}. However, existing transformer-based models can suffer from computational inefficiency when processing large sequences (e.g., text, images, or genomics data) and performing tasks where long-range dependencies are important. Structured state-space sequence models have been proposed as promising alternative\cite{gu2023mamba} to resolve these computational bottlenecks. The core of these sequence models consists of a linear state-space block mapping inputs to  outputs through a latent state, which encompasses information on the past, see Fig.~\ref{fig:mamba-s6-block} \cite{alonso2024state}. The system matrices can be further chosen as input-dependent to propagate or forget information of long sequences. Moreover, other pre- and post-processing (e.g., convolution, activation, or softmax operations) blocks are leveraged. When such structured state-space sequence models are stacked, the overall neural network architecture achieves state-of-the-art performance in various  domains\cite{gu2023mamba}.
\end{example}
}

In these examples, {\tb systems theory has been used to analyze optimization algorithms  or  foundation models that can be abstracted 
as feedback interconnections of few blocks.
There are many other applications} where the full power of a system-theoretic approach can be brought to bear. For example, computing pipelines discussed also in Section~\ref{subsec:architectures} comprise a cascade of algorithms that perform data analysis tasks (e.g. data collection, aggregation, stochastic gradient updates, model deployment). While each of the individual blocks can be understood in isolation, {\tb their interconnection may lead to unexpected behaviors -- a feature that system theorists are well accustomed with \cite{van2000l2,vslijak2007large,hill2022dissipativity}.} An important example are deep-learning-based recommender systems, where current pipelines are brute-force optimized via stochastic gradient descent \cite{ning2022eana} 
and may lead to sub-optimal performance, frequent crashes, or require impractical cluster reorganizations \cite{nagrecha2023intune}.
These last examples leave the realm of algorithms confined to the computing world towards algorithms that interact with their environment. We delve into this topic in Section~\ref{eq:RTAFL} below.

\subsection{Real-Time Algorithms in Feedback Loops}
\label{eq:RTAFL}


Increasingly, algorithms dynamically engage in real-time with real-world scenarios, in closed-loop with physical plants, social and socio-technical networks, data-generating processes, or even other algorithms.
For instance:
\begin{itemize}
\item In physical systems, iterative optimization algorithms are essential for optimization-based control. Such algorithms continuously adapt and respond in real time to changing operating conditions and exogenous disturbances, maintaining optimal performance of the physical plant.

\item Within the sphere of human interaction,  as in social networks, similar algorithms power recommender systems, shaping user experiences based on evolving individual preferences, behavioral patterns, and global trends. The resulting user serves again for training of the algorithms.

\item In generative AI, iterative training algorithms learn from vast, ever-changing datasets and enable the creation of increasingly sophisticated outputs (e.g., realistic images, life-like synthesized voices, accurate predictive text). The generated data eventually serves again for training.

\item In reinforcement learning, {\tb candidate control policies are evaluated in roll-outs. The resulting reward signals are then used to update the  policies, resulting in a feedback loop between policy evaluation and policy improvement.}

\item In the domain of networked systems, agent-based algorithms play a pivotal role in parallel computing, distributed optimization, and finding equilibria of non-cooperative games, facilitating efficient information processing and decision-making across network nodes. 
\end{itemize}

Additionally, algorithms engage in feedback interactions with other algorithms, for example in bi-level programs, where one algorithm's output becomes another algorithm's input, creating a multi-layered decision-making process. Finally, there are also scenarios where algorithms dynamically interact in a fully competitive setting, e.g., in high-speed trading.

\begin{figure}[t]
\centering
\includegraphics[width=.95\columnwidth]{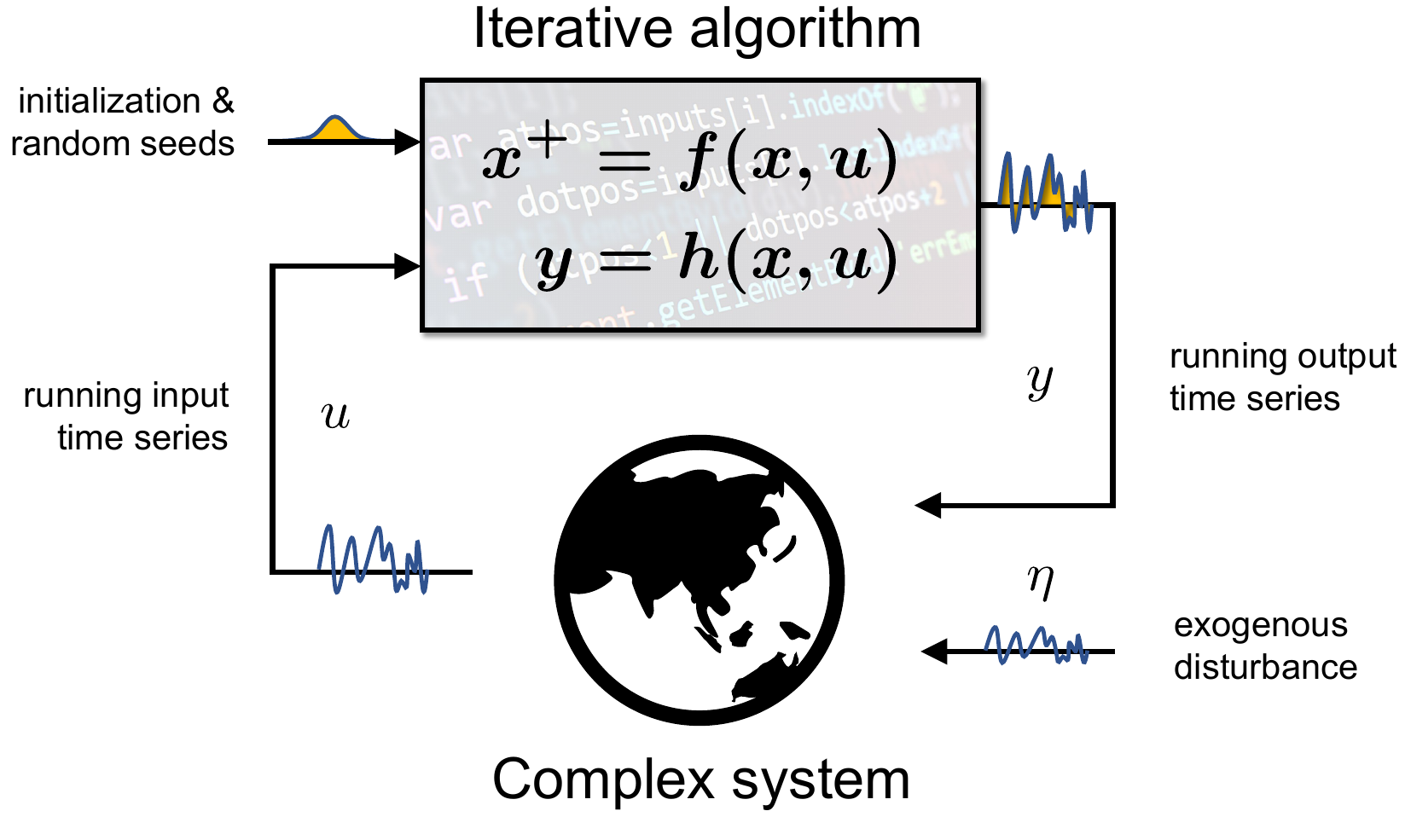}
\vspace{-1em}
\caption{Contemporary iterative algorithms emerging in different research areas and application domains run in real-time and closed-loop with complex (i.e., networked, uncertain, noisy) dynamical systems, such as physical plants, social networks, data-generating processes, or other iterative algorithms. \label{fig:FI}}
\vspace{-2em}
\end{figure}

Viewed through the lens of systems theory, these {\em in vivo} algorithms can be conceptualized as {\em feedback interconnections} between dynamical systems, as illustrated in Figure \ref{fig:FI}.
Here, the ``world'' block can assume different meanings depending on the specific context. For example, in optimization-based control schemes (e.g., MPC \cite{diehl2005real, liao2020time,zeilinger2011real,borrelli2017predictive}, {\tb extremum seeking \cite{ariyur2003real,scheinker2024100}}, feedback optimization \cite{AH-SB-GH-FD:21} and equilibrium seeking \cite{belgioioso2022online}), the ``world'' embodies the physical plant that we would like to operate safely and efficiently. In distributed optimization and network games, it represents the communication protocol, such as dynamic average consensus \cite{zhu2010discrete}, used by the nodes to compute network-level quantities, such as the average gradient in gradient-tracking schemes \cite{nedic2017achieving}, {\tb \cite{carnevale2024unifying}} or the average strategy in aggregative games \cite{belgioioso2020distributed}. 
In reinforcement learning, for example in policy gradient methods like REINFORCE \cite{reinforce}, the world represents the plant on which a policy can be tested to assess the corresponding reward, which is then used in an iterative search procedure; see \cite{Benbrahim1994} for a seminal work where the plant is a robot.
Finally, in hyper-gradient schemes for bi-level programs \cite{grontas2023big}, the ``world'' represents the iterative algorithm used to solve the lower-level optimization problem.

The common challenges in all these scenarios stem from the complex and often real-time feedback interplay between the iterative algorithms and the dynamical systems they are interconnected with. Traditional asymptotic performance metrics for iterative algorithms and  convergence certificates for dynamical systems fall short in such dynamic environments. The focus instead shifts to evaluating input-output properties (e.g., robust stability), understanding error propagation (e.g., via gains), considering transient metrics (e.g.,  running cost or regret), and deriving quantitative certificates on the interconnection of systems at different time scales (e.g., via singular perturbation analysis). These factors demand a nuanced approach to gauge algorithmic performance in online and noisy contexts, where algorithms and systems interact continuously.

Building on this understanding, systems theory emerges as a common {\em lingua franca} and a formidable framework to tackle these challenges. It brings to the fore mathematical tools such as dissipativity, small-gain theory, and singular perturbation methods {\tb cultivated in our community \cite{van2000l2,khalil2002nonlinear,vslijak2007large,hill2022dissipativity}, which are instrumental in examining and ensuring the stability and robustness} of feedback interconnections. 
%
This system-theoretic perspective has been crucial in many success stories across various applications, from control systems to network strategies and data-driven algorithms. We outline some illustrative examples and provide references for a detailed exposition.

%

\begin{example}[Sub-optimal Model Predictive Control]
\label{ex:MPC}
Standard MPC is a feedback law $\kappa(x_t)=\nu_0$ that given the state, $x_t$, at time $t$, returns the first sample of the control input sequence, $\nu$, that solves a multi-stage optimization problem of the form
\begin{align}
\hspace*{-5pt}
S(x_t) = \left\{
\begin{array}{r l l}
\displaystyle
\arg \min_{\xi,\nu} \quad & \sum_{k = 0}^{K-1} \ell_{\text{s}}(\xi_k,\nu_k) + \ell_{t}(\xi_K)\\
\text{s.t.} \quad & \xi_{k+1} = f(\xi_k,\nu_k),\qquad  \forall k \in \mathcal H \\
               & \nu_k \in \mathcal U, \hspace*{6.1em}  \forall k \in \mathcal H \\
               & \xi_0 = x_t,
\end{array}
\right.
\label{eq:POCP}
\end{align}
where $ \ell_{\text{s}}$ and  $\ell_{\text{t}}$ are the stage and terminal costs, respectively, $f(\cdot,\cdot)$ encodes the plant dynamics, $\mathcal U$ are the input constraints, and ${\mathcal H}=\{0, \ldots, T-1\}$. As the problem \eqref{eq:POCP} is parametrized in the initial conditions $x_t$, it needs to be re-solved at every sampling time.
Solving \eqref{eq:POCP} to optimality in real time may not always be feasible for systems with limited computing power, fast sampling rates, network structure, or highly nonlinear dynamics.
Nevertheless, $S(x_t)$ can be approximated by maintaining a running solution estimate $z_t$ and improving it at each sampling instant by using $n$ iterations of an iterative optimization method $\mathcal T$ (e.g., sequential quadratic programming, as in \cite{Diehl2005})  warm-started with the estimate from the previous sampling instant $z_{t-1}$. This gives rise to an algorithm
\begin{align}\label{eq:iterative}
z_t = \mathcal T^n(z_{t-1},x_t) 
\end{align}
From a system-theoretic perspective, this sub-optimal MPC law can be conceptualized as a feedback interconnection between two dynamical systems, the optimization algorithm and the plant, as illustrated in Figure \ref{fig:soMPC}.
The stability and robustness of this cyber-physical interconnection can be established using a celebrated system-theoretic result: the \textit{small-gain theorem}.
Remarkably, this system-theoretic analysis reveals that there is an (a-priori) fixed number of iterations $n$ of the algorithm $\mathcal T$ for which sub-optimal MPC retains similar stability and robustness properties to its nominal counterpart\,\cite{liao2020time}. {\tb Intuitively, this can be also regarded as a time-scale separation condition requiring the optimization dynamics \eqref{eq:iterative} to be sufficiently faster than the plant dynamics.}%
\begin{figure}[h]
\centering
\includegraphics[width=.95\columnwidth]{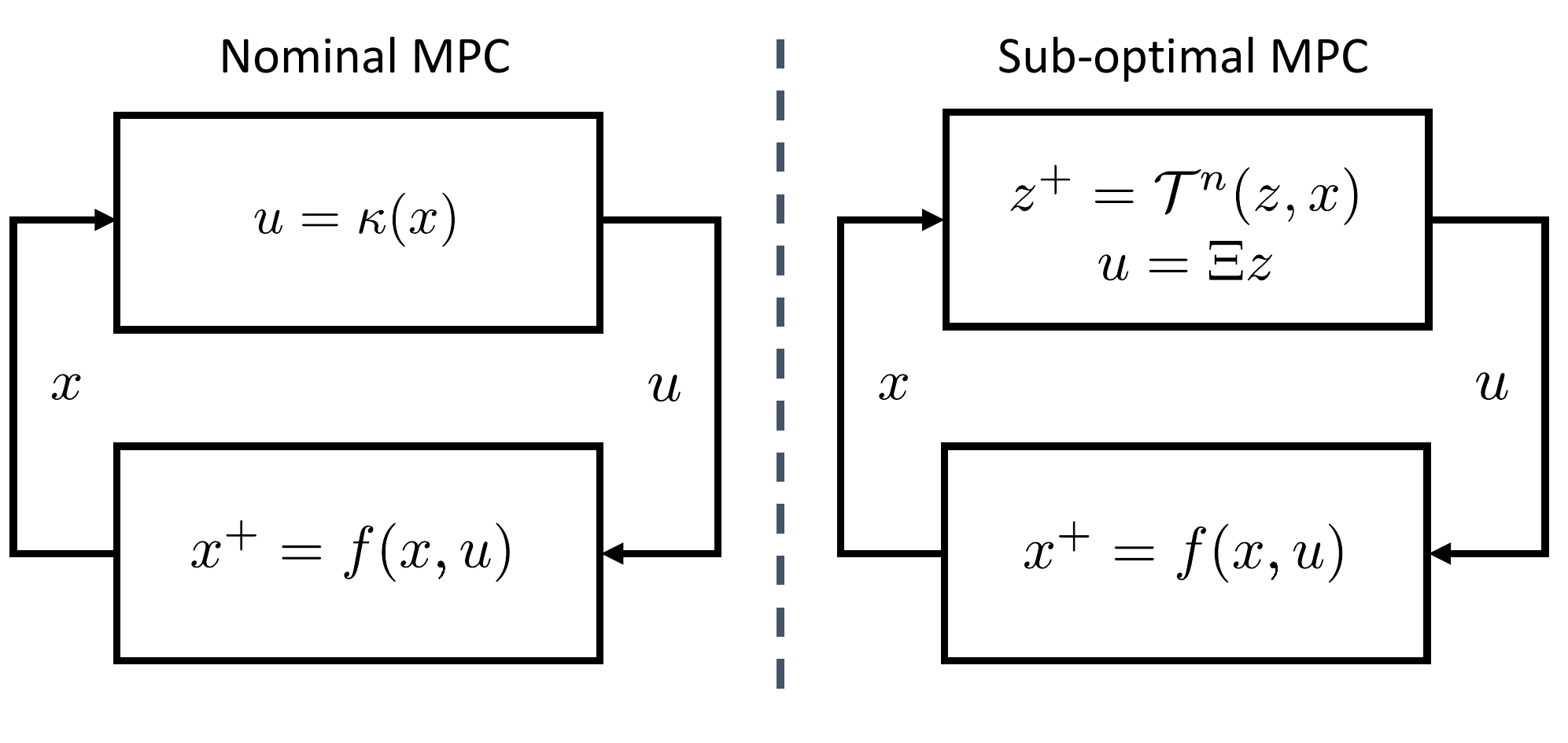}
\vspace{-.5em}
\caption{\cite[Figure 1]{liao2020time} Left: Nominal MPC as a static feedback law $\kappa(x)$. Right: Sub-optimal MPC as a feedback interconnection between a physical plant and an optimization algorithm with a solution estimate $z$, as its internal state, dynamics defined by its iteration $\mathcal T(z,x)$, and an output matrix $\Xi$ selecting the first sample of the control input sequence $z$.\label{fig:soMPC}}
\vspace{-1em}
\end{figure}
\end{example}

\begin{example}[Distributed Optimization via Gradient Tracking]
\label{ex:GTrack}
In this example, extrapolated from \cite{nedic2017achieving}, we consider the following distributed convex optimization problem
\begin{align}
\label{eq:DOP}
\min_{x} \; \frac{1}{n} \sum\nolimits_{i =1}^n f_i(x),
\end{align}
where each function $f_i$ is held private by agent $i$. The agents are connected via a communication network and want to collaboratively solve \eqref{eq:DOP}, while exchanging information with neighboring agents. A conceptual iteration solving \eqref{eq:DOP} is
\begin{align}
\label{eq:CI}
\boldsymbol x^{k+1} = W \boldsymbol x^k - \alpha \frac{1}{n} \mathbf 1 \mathbf 1^\top  \nabla f(\boldsymbol x^k),
\end{align}
where $\boldsymbol x = {[x_1^\top \; \ldots \; x_n^\top]}^\top$ is the stacked vector of local solution estimates $x_i$, $W$ is a doubly-stochastic mixing matrix that is consistent with the network structure, $\alpha$ is a step size, and $\nabla f(\boldsymbol x) = [\nabla f_1(x_1)^\top \; \ldots \; \nabla f_n (x_n)^\top]^\top$  collects the local gradients.
%
%
Unfortunately, \eqref{eq:CI} is not amenable to a distributed implementation, as it
requires a central node to provide the average of the gradients $\frac{1}{n} \mathbf 1 \mathbf 1^\top  \nabla f(\boldsymbol x)$. 
Nevertheless, a distributed approximation of  \eqref{eq:CI} can be obtained by replacing the average gradient with surrogate variables $\boldsymbol s$, updated as
\begin{align}
\boldsymbol s^{k+1} &= W \boldsymbol s^k + \big(\nabla f(\boldsymbol x^{k+1}) - \nabla f(\boldsymbol x^{k})\big),
\label{eq:DAC}
\end{align}
and initialized as $\boldsymbol s^0 =\nabla f(\boldsymbol x^{0})$. The update \eqref{eq:DAC} is an instance of dynamic average consensus, a technique originally introduced in \cite{zhu2010discrete} for dynamically tracking the average state of a multi-agent system. {\tb After the change of variables $\boldsymbol z = \boldsymbol s - \nabla f(\boldsymbol x) $, the resulting distributed iteration reads as
\begin{subequations} 
\begin{align}
\label{eq:OPTDYNAM}
\boldsymbol x^{k+1} &= W \boldsymbol x^k - \alpha ( \nabla f(\boldsymbol x^k)  + \boldsymbol z^k),\\
\label{eq:CONDYNAM}
\boldsymbol z^{k+1} &= W \boldsymbol z^k - (I-W)  \nabla f(\boldsymbol x^k),
\end{align}
\end{subequations}
with $\boldsymbol z^0 =0$.}
This can be seen as a feedback interconnection between the optimization algorithm and the dynamic average consensus protocol, as illustrated in Fig.~\ref{fig:DIG}. {\tb If the consensus dynamics \eqref{eq:CONDYNAM} are sufficiently fast, then the optimization dynamics \eqref{eq:OPTDYNAM} track the nominal algorithm in \eqref{eq:CI}. This intuition can be exploited with formal small-gain \cite{nedic2017achieving} and (similarly) singular perturbation \cite{carnevale2024unifying} arguments that yield convergence certificates, including bounds on the step size $\alpha$. }
\begin{figure}[h]
\vspace{-1.5em}
\centering
\includegraphics[width=.55\columnwidth]{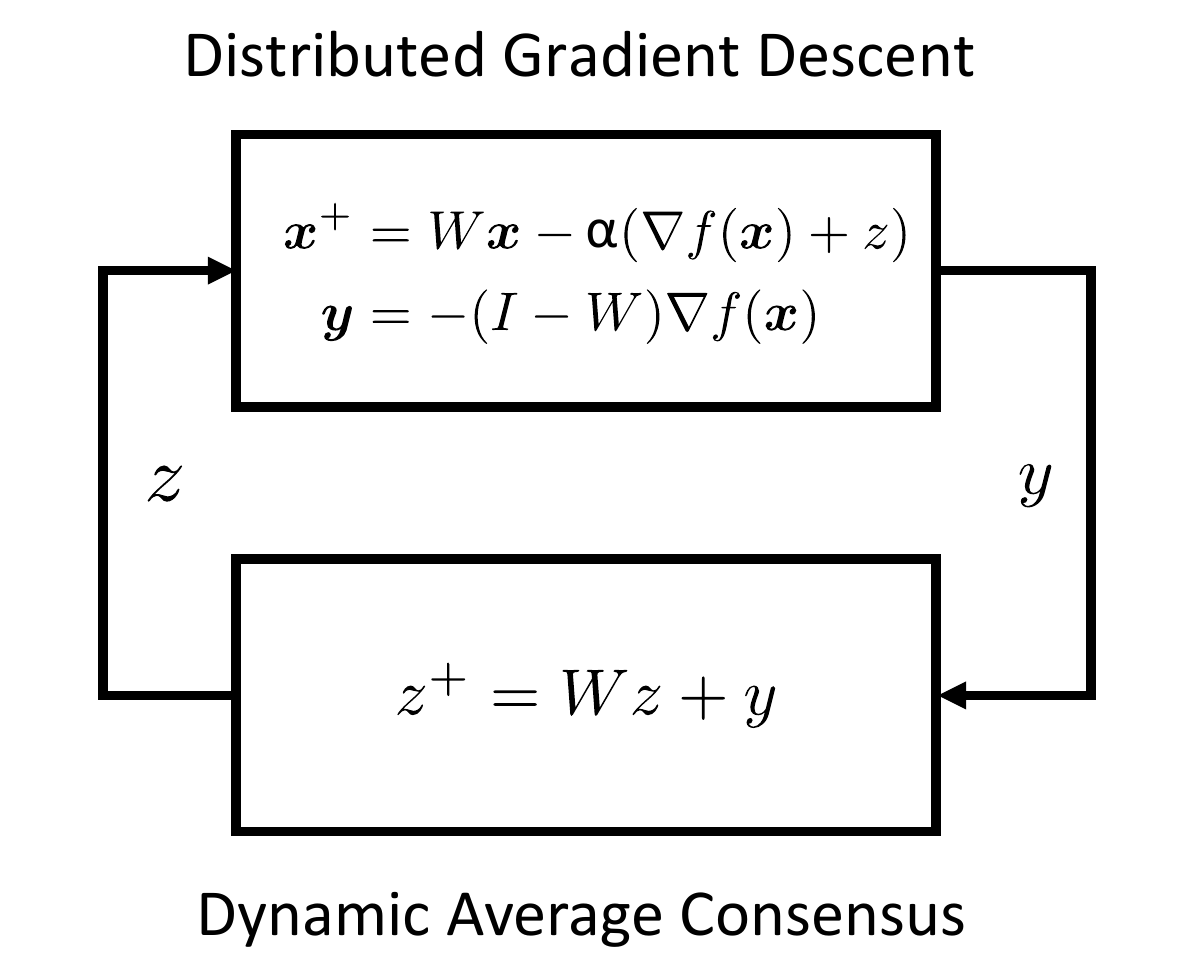}
\vspace{-1em}
\caption{A system-theoretic interpretation of the {\em gradient-tracking} algorithm \cite{nedic2017achieving} as a feedback interconnection between the gradient descent dynamics and a dynamic average consensus protocol \cite{zhu2010discrete}. The dynamics are obtained from \eqref{eq:CI}--\eqref{eq:DAC} with the change of variables $\boldsymbol z = \boldsymbol s - \nabla f(\boldsymbol  x)$, see \cite{bin2019system}.\label{fig:DIG}}
\vspace{-1em}
\end{figure}
\end{example}

\begin{example}[Online Feedback Optimization]
\label{ex:OFO}
Consider the problem of ``efficiently" operating a physical plant described by the state-space system
\begin{subequations}
 \label{eq:ctime-dynamics}
 \begin{align}
\dot{x} &= f(x,u,w),\\
y &= g(x,w),
\end{align}
\end{subequations}
where $x$ is the state, $y$ is the output, $u$ is the control input, and $w$ is an exogenous disturbance. We assume that \eqref{eq:ctime-dynamics} is stable and admits a steady-state input-output mapping $ h(u,w)$.
%
The control objective is to steer \eqref{eq:ctime-dynamics} to an economic steady state, implicitly described via the optimization problem
\begin{subequations} 
\label{eq:OFO_ss}
\begin{align} 
\min_{u,s} \quad & \Phi(u,s)\\
\textrm{s.t.} \quad & s = h(u,w),\\
                   & u \in \mathcal{U}, \; s \in \mathcal{Y}
\end{align}
\end{subequations}
where $\Phi(u,s)$ evaluates the steady-state performance of the plant, $\mathcal U$ and $\mathcal{Y}$ are the sets of admissible inputs and outputs.

If the input-output model $h(\cdot)$ is perfectly known, and the exogenous signal $w$ is perfectly predicted, then \eqref{eq:OFO_ss} can be solved offline, and the solution implemented on the system \eqref{eq:ctime-dynamics}. This {feedforward} approach, {\tb represented in the left panel of Figure~\ref{fig:OFO},} lacks robustness to model uncertainty and variations in the exogenous disturbance.
Online Feedback Optimization \cite{AH-SB-GH-FD:21} is an alternative approach to this steady-state regulation problem driven by the system-theoretic intuition that feedback introduces robustness.
The core idea is {\tb to steer the dynamical system \eqref{eq:ctime-dynamics} to the solution of the} optimization \eqref{eq:OFO_ss} by using any (appropriate) iterative algorithm  of the form
\begin{align}
\label{eq:OFO}
u^{k+1} = \mathcal T(u^k,y^{k}).
\end{align}
{\tb Crucially, in \eqref{eq:OFO} the evaluation of the input-output model $h(u,w)$ is replaced by online measurements of the output $y$ obtained directly from the plant \eqref{eq:ctime-dynamics}.}
This closed-loop optimization procedure not only improves robustness against time-varying disturbances but also reduces model dependence and computational effort, as the evaluation of $h(u,w)$ is outsourced to the physics.
{\tb When \eqref{eq:OFO} takes the form of a stochastic gradient descent, the scheme resembles classical extremum seeking \cite{ariyur2003real,scheinker2024100}.}
From a system-theoretic perspective, Online Feedback Optimization is a feedback interconnection of two dynamical systems: the iterative algorithm and the physical plant, as shown in the right panel of Figure \ref{fig:OFO}. 
This system-theoretic perspective applies more generally to \emph{Feedback Equilibrium Seeking} problem of steering a dynamical system so that it tracks the solutions of time-varying generalized equations; these can be local minimizers of nonlinear programs as in \eqref{eq:OFO_ss} but also {\tb Nash equilibria} of non-cooperative games.
Input-to-state stability and small-gain theory prove to be powerful tools to certify the stability and robustness of this cyber-physical interconnection \cite{belgioioso2022online}.

\begin{figure}[h]
\centering
\includegraphics[width=\columnwidth]{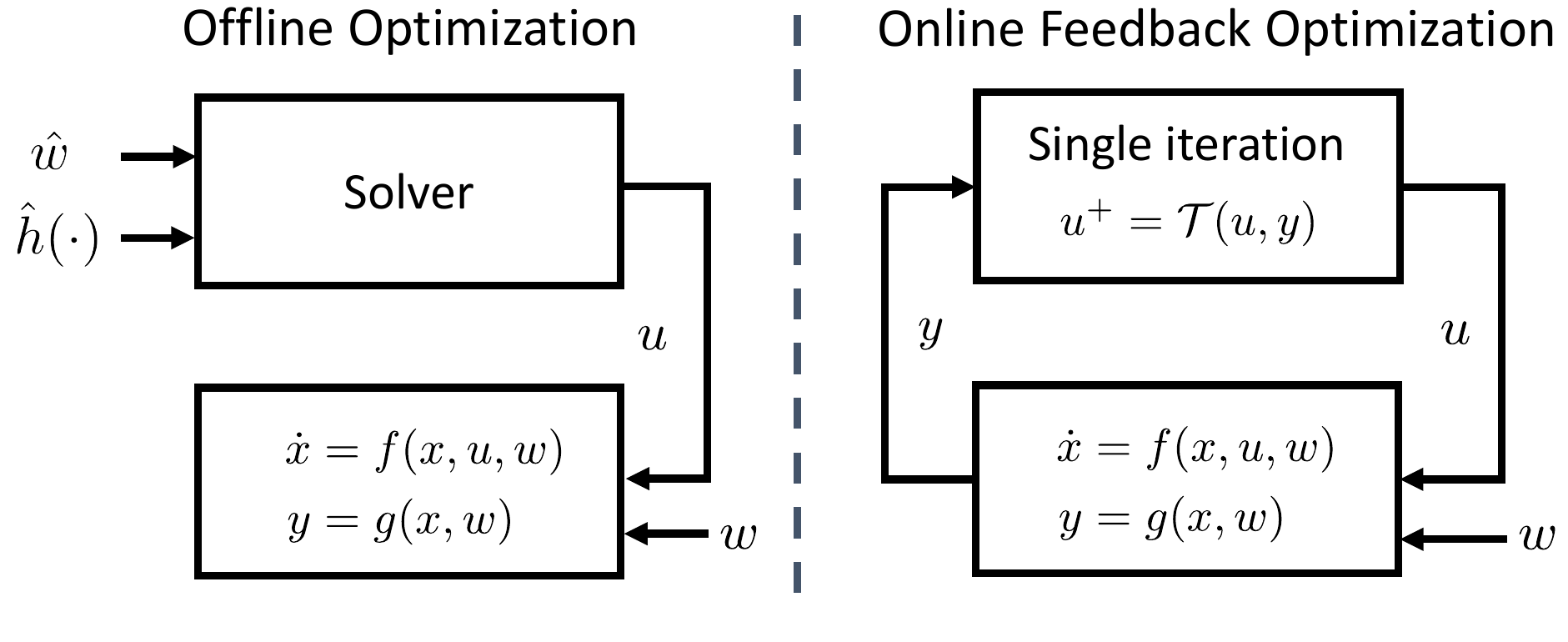}
\vspace{-1.5em}
\caption{Left: In offline optimization, a static model $\hat h$ of the physical plant \eqref{eq:ctime-dynamics} and a forecast $\hat w$ of the exogenous disturbance are used within an appropriate solver to find a solution to \eqref{eq:OFO_ss} offline; then, the outcome of the offline optimization routine is fed to the physical plant. Right: In online feedback optimization, measurements $y$ obtained by sampling the physical plant \eqref{eq:ctime-dynamics} are integrated within an iterative algorithm $\mathcal T(\cdot)$ that solves \eqref{eq:OFO_ss}, creating a cyber-physical feedback interconnection robust to exogenous disturbances $w$. \label{fig:OFO}}
\vspace{-1.5em}
\end{figure}
\end{example}

These examples demonstrate the profound impact that a system-theoretic perspective can provide in enhancing our analysis and design of real-time algorithms in feedback loops. 
The small-gain analysis in Example~\ref{ex:MPC} paves the way for applying MPC \cite{borrelli2017predictive} -- the go-to methodology for constrained control -- to a broader array of safety-critical systems that cannot benefit from standard MPC implementations, due to limited computing power, fast sampling rates \cite{liao2020time,zeilinger2011real}, or large network structures \cite{belgioioso2023stability}.
Similarly, the intuition of integrating dynamic average consensus into (centralized) algorithms, as showcased in Example~\ref{ex:GTrack}, has been instrumental in the development of many novel distributed algorithms for multi-agent optimization and noncooperative games over large networks, {\tb (see \cite{carnevale2024unifying} and references therein)}.
Finally, the online feedback optimization approach in Example~\ref{ex:OFO}, by allowing the optimal operation of complex physical systems despite the lack of precise model information, has provided a novel solution to several engineering problems, especially in the real-time operation of electrical power systems \cite{simonetto2020time,AH-SB-GH-FD:21,OrtmannAEW,kroposki2020autonomous,dall2016optimal}.

To conclude this section, we highlight some areas where the system-theoretic perspective, while not fully harnessed or widely accepted, has the potential to provide new insights. 
One exciting possibility is using a control-theoretic approach to study the phenomenon of performative prediction, where the decision of an algorithm (typically a predictive one) affects the environment and, therefore, the data that is fed back into the same algorithm. {\tb The predominant paradigm of statistical learning addresses data-generating distributions unaffected by algorithmic decisions, excluding the phenomenon of endogenous distribution shift, i.e., performativity \cite{hardt2023performative}.}
The inability to understand and control this phenomenon is particularly concerning in socio-technical systems, where it is closely connected to the emergence of bias and unfairness \cite{Pagan2023}.
For example, machine learning predictions used in loan approvals impact future default rates that, in turn, affect the data that is used for training the machine learning models. Similarly, we are currently witnessing the deployment of large language models already significantly impacting content on the internet, and these will be again the basis for training new language models. Systems theory has the potential to provide a framework for understanding interconnections between computation, predictions, and interactions with the real world.


Another inspiring direction comes from the field of data-driven predictive control \cite{IM-FD:22-tutorial,Berberich2021}.
This control methodology has recently gained research momentum thanks to its ability to provide predictive control actions directly from raw data, thus bypassing the often expensive system identification step \cite{FD:23-editorial2}. However, this comes at the price of an increased complexity of online computations. This curse of dimensionality often renders data-driven predictive control inapplicable to large-scale and/or highly nonlinear systems, where, ironically, model-free control would be most beneficial. Drawing parallels with the MPC strategy from Example~\ref{ex:MPC}, a potential solution may lie in adopting a suboptimal version coupled with a similar small-gain analysis. This approach holds promise for mitigating the current computational challenges, thereby unlocking the potential of data-driven predictive control for broader applications in complex environments, such as infrastructure networks.

\subsection{Decision-Making Architectures}\label{subsec:architectures}

Control engineers are inherently inclined to abstract and conceptualize in block diagrams, imposing structure, modularization, and managing complexity. This leaves them uniquely positioned to analyze and design architectures emerging in algorithmic pipelines, extensive machine learning systems, or autonomous decision-making.
For instance, we are accustomed to addressing queries like determining the demarcation of system boundaries, identifying the interacting components, assessing the extent of tolerable uncertainty, ensuring desired certifications under interconnected systems, and similar considerations. 
These questions are increasingly relevant in large algorithmic pipelines
and even within individual  components\footnote{For example, DARPA's {\em Data Driven Discovery of Models} program is concerned with the automation of machine learning pipelines both in terms of their architecture and components: \url{https://datadrivendiscovery.org}.}.%

\begin{figure}[tb]
	\centering{
	\includegraphics[width=.85\columnwidth]{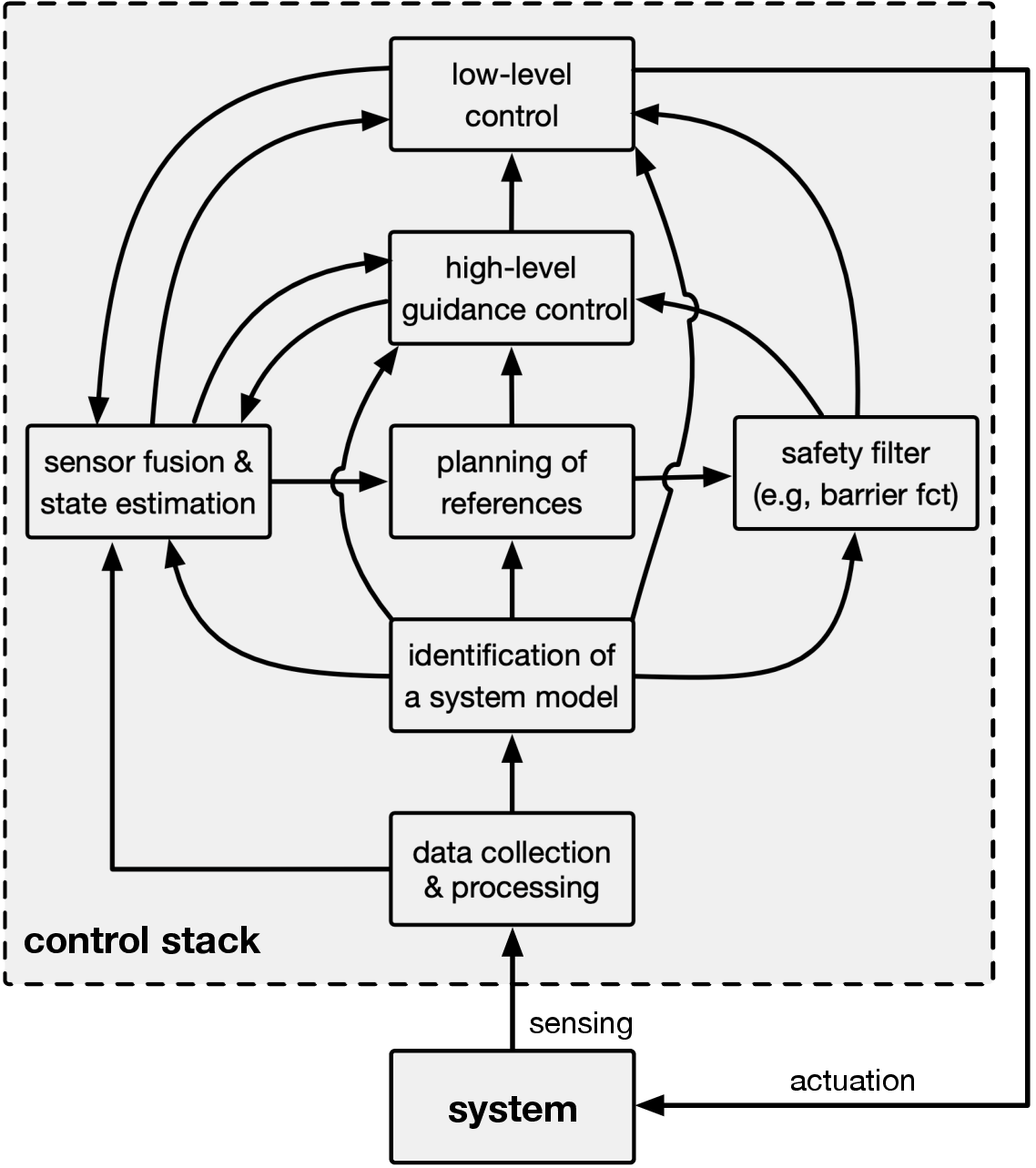}
		\vspace{-.5em}
	\caption{
 The block diagram illustrates a prototypical layered decision-making architecture -- also termed control stack -- in an automatic control system.
 }
	\label{Fig: layered stack}
	}
	\vspace{-1.5em}
\end{figure}
Rather than arguing in general terms, let us focus on interconnected algorithms and architectures in control systems, which is a running thread in the IEEE Control System Society 2030 roadmap \cite{css-roadmap}. 
Fig. \ref{Fig: layered stack} depicts a prototypical layered decision-making architecture --- sometimes also referred to as a {\em control stack} or {\em autonomy stack} --- implementing an automatic control system from sensing, over identification, estimation, planning, and so forth, to the ultimate actuation of the system, and finally closing the loop back to sensing. We highlight that every block in Fig. \ref{Fig: layered stack} relies on an iterative numerical algorithm (unless it is explicit such as a PI  control) prone to propagating uncertainty or errors especially when implemented in real-time or resource-constrained platforms.%

While the architecture in Fig. \ref{Fig: layered stack} is prototypical and often encountered with minor deviations, we remark that many competing architectures have been put forward that either merge or side-step some of the elements. Examples include:
\begin{itemize}

\item {\tb (Economic) MPC \cite{faulwasser2018economic,borrelli2017predictive} merges planning and control, and it can also naturally include moving horizon estimation as well as the safety filter tasks -- the latter of which is often separately handled by barrier functions \cite{wabersich2023data}.} 

\item Direct data-driven control methods side-step the  identification of a model or condense the entire control stack into a single end-to-end learning problem of the form
\begin{align}
    \text{minimize}_{u,y} &\;\text{ objective}(u,y)
        \label{eq:e2e}
    \\
    \text{subject to}& \;\; (u,y) \text{ compatible with recorded data}
    \nonumber
\end{align}
where $(u,y)$ denote the inputs and outputs.\footnote{We are deliberately vague concerning the ``compatible'' quantifier in \eqref{eq:e2e}. It may literally correspond to minimizing a residual between $(u,y)$ and observed data, as in imitation learning \cite{osa2018algorithmic} or data-enabled predictive control \cite{IM-FD:22-tutorial}.}
The conventional, layered architecture for solving \eqref{eq:e2e}, starting with system identification and possibly including other elements from Fig. \ref{Fig: layered stack}, restricts the search space of  control policies in \eqref{eq:e2e} and facilitates the solution by imposing structure, injecting side-information, filtering of noise, and so on. It may also, however, induce erroneous bias and suboptimality; see \cite{FD:23-editorial2} for a discussion.

\end{itemize}

Some core foundations of our discipline concern separation principles. Despite looking intrinsically coupled, different subroutines in Fig. \ref{Fig: layered stack} can be designed and executed separately; a typical example is LQR and Kalman filtering. Conversely, some subroutines that look inherently sequential  may benefit from bi-directional interactions, e.g., identification for control recognizes that the ``downstream'' control task should bias the ``upstream'' identification task \cite{hjalmarsson2005experiment}. Similar considerations can be found in the machine learning and operations research communities under different names \cite{mandi2023decision,sadana2023survey}.

A relevant historical example  of interconnected algorithms akin to a control stack is internet congestion control, where the system-theoretic approach brought structure, interpretability, and improvements to the architecture \cite{low2002internet,KELLY2003159,chiang2007layering}. Communication over the internet is governed by transmission control protocols which, in their simplest form, are based on the following mechanism: A source probes the network by increasing its transmission rate until congestion is detected, at which point the rate is decreased.
Crucially, the lower-level protocols that control the transmission of the data stream over the physical media need to be able to provide upstream information about congestion (which can be detected through delays, dropped packets, or explicit packet flags set by intermediate routers).
This mechanism can be reverse-engineered as a saddle-point seeking dynamical system (similar to \eqref{eq:saddle_point_dynamics}) that converges to the solution of a network utility maximization problem. This system-theoretic perspective on internet congestion control provided important insights into optimality, stability, and robustness with respect to network capacity and delays.

The above examples provide plenty of arguments in favor of one or the other architecture, and a truly {\em grand challenge} concerns the layering and decomposition of the control stack \cite{css-roadmap,matni2024towards}. We believe that a system-theoretic perspective, namely viewing an algorithm stack as interconnected systems, provides a promising angle on how to draw and analyze the block diagram, e.g., concerning propagation of uncertainty, (sub)optimality of architectures, or decomposition into subtasks, among others. Recent notable algorithmic perspectives on the control stack concern the separation and coupling of planning and tracking \cite{srikanthan2023augmented} or the feedback interconnection of recursive least squares and policy iteration \cite{song2024role}. As  a last thought, the whole composition of algorithms making up the control stack becomes itself an algorithmic system which can be analyzed as such concerning its input/output properties. 


{\tb We did not touch upon several ``digital'' algorithmic aspects related to Fig. \ref{Fig: layered stack}: e.g., different levels of the control stack typically run at different rates, employ models of different fidelity, and are formulated in either discrete or continuous time; see \cite{matni2024towards}. Further, Fig. \ref{Fig: layered stack} does not display the higher {\em logic} algorithmic layer whose specifications dictate the choices of optimization criteria and constraints in individual blocks in Fig. \ref{Fig: layered stack}. This logic layer comes with its own architectural challenges and solutions \cite{alur1999reactive,keutzer2000system} that need to be integrated.}

Undoubtedly, within our discipline, there is a need to comprehend and devise an effective approach to architectural challenges, such as layering and modularity. We believe that this discussion and the ensuing methodologies will be beneficial not just within control but also for broader applications in interconnected algorithms and pipelines.

%
%
%
%
%


\section{Challenges that Appear Approachable from the System-Theoretic Perspective} 
\label{sec:challenges}

We now examine the potential of the system-theoretic perspective to offer new solutions to contemporary challenges of {\tb algorithms operating {\em in silico} or {\em in vivo}.} 
Our focus is on introducing a novel perspective that opens up fresh avenues for tackling a range of complex problems, aiming to surpass the state of the art. Going beyond the case studies discussed in the earlier sections (where systems theory has already achieved specific successes, e.g. in optimization algorithm design \cite{van2017fastest}), we identify overarching problems and highlight how systems theory can lead to significant insights and advancements.
%

{\tb
\paragraph{Modeling} 
 Adopting a system-theoretic approach enables valuable abstractions, 
 particularly in modeling the closed-loop interactions between machine learning algorithms and their (physical, social, or biological) environments: how do machine learning predictions impact their environments, and how do streams of 
 closed-loop data affect  learning outcomes? Often times it is unclear how to draw the block diagram and identify  feedback loops.
  To our knowledge, such questions are only beginning to be understood in static settings \cite{hardt2023performative}.
Further, often probability is a natural lingua franca to model the world in presence of uncertainties, priors, or when describing the collective behavior of a large population. One use case, where all of the above come together, is the feedback interplay of recommender systems and a user population \cite{lanzetti2023impact}. To model these population-level interactions, we envision fruitful research avenues surrounding algorithms interacting with an environment modeled intrinsically on probability spaces, which is attracting growing interest in the control community~\cite{chen2021optimal}.

\paragraph{Breaking the hierarchy in  interconnections} Several optimization schemes and optimization-based control methods rely on the interconnection between two (or more) iterative algorithms and their (physical) environment.
A pivotal yet restrictive assumption (often inherent in these interconnected dynamical systems) is their operation within distinct time scales (e.g., as illustrated in Example~\ref{ex:MPC}). On the algorithmic side, this often manifests itself as nested {\tt for} loops; when interconnected with a physical plant, this corresponds to nearly instantaneous dynamics. 
%
System-theoretic tools -- such as  singular perturbation or small-gain methods, feedforward, backstepping, and dissipativity certificates -- present promising avenues to dismantle this hierarchical structure; see \cite{picallo2022sensitivity} for applications in optimization. 
Extrapolating from the examples in Section~\ref{sec:examples}, concrete questions concern co-design (of communication and optimization protocols) in distributed optimization, (of nested loops) in multi-level optimization, or (of stabilizing controllers and optimization algorithms) in feedback-based optimization to avoid time-scale separation.

\paragraph{Architecture \& composition aspects} The system-theoretic viewpoint is suited to understand feedback effects and compositions of algorithms in learning systems. There are various channels (e.g., data collection, sampling mechanism, and feature extraction \cite{Pagan2023}) through which the downstream decision influences the upstream problem. Quantifying, altering, and utilizing such feedback effects are promising topics. Moreover, compositions of algorithms may lead to cascades (e.g., machine learning pipelines), parallelization (e.g., federated learning), feedback interconnections (e.g., bilevel optimization), or other types of layering as in the control stack in Fig.~\ref{Fig: layered stack}. Key challenges are analyzing and designing composition structures to improve the overall performance, including convergence rates, robustness, and generalization.
%
Further, system-theoretic analysis can help to characterize how errors and uncertainties propagate through composition of algorithms, as shown early on for MPC \cite{diehl2005real} or recursive least squares \cite{ljung1985error}. Leveraging familiar system-theoretic small-gain arguments, we may expect that if the propagation ratio of each block is less than a threshold, the global pipeline will not suffer from divergence due to cumulative errors or uncertainties. However, a small gain may also intricately imply a slow convergence rate. It is important to explore suppression mechanisms and interconnection certificates that strike a balance between robustness and fast convergence.

\paragraph{Beyond analogue block diagrams} Our presentation was mostly idealized and aside from brief disclaimers neglected ``digital'' challenges in the implementation of algorithms like multi-threading, asynchronous communication and computation, quantization, mixed discrete and continuous time, multiple sampling rates, and so on. Further, the examples put forward in this paper all consider numerical algorithms operating on real numbers. However, some algorithms cannot immediately be posed in this framework, e.g., the simplex method, or they operate over discrete (or non-numerical spaces) such as combinatorial optimization, sorting of lists, graph search, and so on. While some of these problems can in principle be solved by continuous dynamical systems \cite{brockett1991dynamical}, further research efforts are necessary to determine whether systems theory offers an effective paradigm for solving such discrete problems. Last, system-theoretic approaches towards algorithms can and should incorporate various digitalization and computer engineering issues and, of course, also integrate computer science elements, such as, e.g., temporal logic.

\paragraph{Realization \& synthesis} In the vast field of optimization, a zoo of algorithms has emerged many of which are known to be {\em equivalent} (with a precise meaning) upon (nonlinear) changes of coordinates, dualization, or embedding into higher-dimensional spaces. In the language of systems theory, such algorithms are different realizations of the same behavior as seen from inputs and outputs. Such a perspective can be leveraged to discriminate between algorithms \cite{zhao2021automatic}, identifying salient features \cite{redman2023equivalent}, or for a data-driven analysis \cite{AP-CJ-HVW:24}. Extrapolating from the  developments of systems theory then leads to quests for normal forms, minimality, I/O analysis and synthesis, and so on. On the last point: most examples in this manuscript concern the analysis of algorithms. Clearly, systems theory also offers a plethora of methods to {\em design} novel algorithms, as demonstrated successfully for a class of optimization algorithms \cite{scherer2023optimization,sundararajan2020analysis,lessard2022analysis}. In a similar spirit, we believe that the synthesis of robust and fast algorithms or their (non-hierarchical) interconnections is one of the most impactful contributions that the control community can offer.
}

\section{Conclusions}
\label{sec:conc}

In this opinion paper, we advocated a view of numerical algorithms as open dynamical systems interacting with their environment. We presented a suite of contemporary and historical examples where this perspective proves useful and listed several promising research directions going forward. 

We envision that {\em systems theory of algorithms} is a topic to be embraced by the systems and control community, as we are uniquely positioned for addressing the related challenges. Our ability to abstract and understand feedback loops offers a unique and powerful perspective for algorithm study and design. Moreover, as modern algorithms increasingly interact with various entities in real-time environments, our deep knowledge of feedback systems positions us to play a pivotal role in this evolving digital and algorithmic landscape. To make a lasting impact, however, we need to broaden our arsenal by incorporating tools not widely used in our field.

We hope that our opinion paper stimulates further research and applications surrounding the systems theory of algorithms.

\section*{Acknowledgements}

The authors wish to thank {\tb the editors, anonymous reviewers, their colleagues (in particular, M. Bianchi, J. Eising, A. Padoan, and D. Liao-McPherson), collaborators, and speakers of the NCCR Symposium "Systems Theory of Algorithms" \cite{nccrsymposium}}
for insightful discussions that led to this paper.

\balance
\bibliographystyle{IEEEtran}
\bibliography{references}

\end{document}